\documentclass[12pt]{article}
\usepackage{latexsym}
\usepackage{amsthm}
\usepackage{amssymb}
\usepackage{amsmath}
\usepackage{amsfonts}
\usepackage{mathrsfs}
\usepackage[all]{xy}
\usepackage{texdraw}
\usepackage{graphicx}
\usepackage{psfrag}
\usepackage{makeidx}

\newtheorem{Theory}{Theory}[section] %Counter for all types of Theorems
\newtheorem{theorem}[Theory]{Theorem}
\newtheorem{lemma}[Theory]{Lemma}
\newtheorem{technicalLemma}[Theory]{Technical Lemma}
\newtheorem{corollary}[Theory]{Corollary}
\newtheorem{proposition}[Theory]{Proposition}
\newtheorem{fact}{Fact}  %level of Exercise
\newtheorem{remark}[Theory]{Remark} %trivial but worth noticing
\newtheorem{question}{Question} %Open Question in my mind or in Theory
\newtheorem{conjecture}[question]{Conjecture}%Statement of suspicion

\newtheorem{Ntn}{Description} %Counter for all types of Notation/Definition
\newtheorem{Dn}[Ntn]{Definition}

\newcommand{\be}{\begin{enumerate}}
\newcommand{\ee}{\end{enumerate}}
\newcommand{\bq}{\begin{question}}
\newcommand{\eq}{\end{question}}
\newcommand{\bcj}{\begin{conjecture}}
\newcommand{\ecj}{\end{conjecture}}
\newcommand{\bc}{\begin{corollary}}
\newcommand{\ec}{\end{corollary}}
\newcommand{\bl}{\begin{lemma}}
\newcommand{\el}{\end{lemma}}
\newcommand{\btl}{\begin{technicalLemma}}
\newcommand{\etl}{\end{technicalLemma}}
\newcommand{\bt}{\begin{theorem}}
\newcommand{\et}{\end{theorem}}
\newcommand{\bp}{\begin{proposition}}
\newcommand{\ep}{\end{proposition}}
\newcommand{\bft}{\begin{fact}}
\newcommand{\eft}{\end{fact}}
\newcommand{\brk}{\begin{remark}}
\newcommand{\erk}{\end{remark}}
\newcommand{\bd}{\begin{Dn}}
\newcommand{\ed}{\end{Dn}}

\newcommand{\ploi}{PL_o(I) }

\newcommand{\Z}{ \mathbf Z }
 
\newcommand{\N}{\mathbf N }
 
\newcommand{\R}{\mathbf R }

\newcommand{\pls}[1]{PL(S^{#1})}
\newcommand{\pli}[1]{PL(I^{#1})}
\newcommand{\Rot}[1]{\mathop{\textrm{Rot}}(#1)}

\author{Collin Bleak}

\begin{document}

\centerline{{\Large{\bf Some questions about the dimension of a group
action.}}}
\vspace{.1 in}

\centerline{Collin Bleak} \centerline{Department of Mathematics,
Cornell University, Ithaca, NY 14853-4201 USA}
\vspace{.3 in}

\centerline{{\bf Abstract}}
\vspace{.1 in} 

We discuss three families of groups, $ZW_n$, $PL(I^n)$, and $PL(S^n)$
(the last two being families of groups of piecewise-linear
homeomorphisms of standard $n$-dimensional spaces).  We show that for
positive $n\in\N$, $ZW_n$ embeds in $PL(I^n)$ which embeds in $PL(S^n)$.
In another direction, $ZW_2$ fails to embed in $PL(I^1)$ by a result
in \cite{bpgsc}.  Here, we extend that result to show that $ZW_2$ also
fails to embed in $PL(S^1)$.  The nature of the proofs of these
non-embedding results leads us to ask if there are corresponding
non-embedding results in higher dimensions.

\tableofcontents 
\newpage
\section{Introduction}

We use a non-embedding result in \cite{bpgsc} to show that the group
$\Z\wr(\Z\times\Z)$ does not embed in $\pls{1}$, the group of
piecewise-linear homeomorphisms of $S^1$ (under the operation of
composition).  This, together with a family of easy embedding results
for a family $ZW_n$ of related groups leads us to ask some natural
questions about the parallel results in higher dimensions.

The author would like to thank Ken Brown for asking whether the
non-embedding of $\Z\wr(\Z\times\Z)$ into $PL_o(I)$ (and therefore
into R. J. Thompson's group $F$) could be extended to a non-embedding
of $\Z\wr(\Z\times\Z)$ into Thompson's group $T$, a subgroup of
$\pls{1}$.  The author is also indebted to Martin Kassabov and
Francesco Matucci for simplyfing the arguments below in accord with
our on-going joint work (see \cite{bkmss}).

\subsection{Families of groups of piecewise-linear self-homeomorphisms}

For each positive integer $n$, let $I_p^n$ and $S_p^n$ represent
polyhedra in $\R^n$ and $\R^{n+1}$ which are homeomorphic to $I^n$ and
$S^n$ respectively, where we further require that $\R^{n+1}$ admits a
reflective plane of symmetry for the polyhedron $S_p^n$, so that
$S_p^n$ is realized as the union of two piecewise-linear homeomorphic
images of $I_p^n$, identified together on their common boundary, where
one image is the reflection of the other across the plane of symmetry,
excepting their boundaries, which lie in the plane of symmetry itself.
We will refer to $\mathrm{Homeo}(I^n)$ and $\mathrm{Homeo}(S^n)$ as
the full groups of self-homeomorphisms of the polyhedra $I^n_p$ and
$S^n_p$ respectively, under composition, and we will refer to
$\pli{n}$ and $\pls{n}$ as the subgroups of these full homeomorphism
groups consisting of the piecewise-linear self-homeomorphisms.

\subsection{Standard restricted wreath products, and another family of groups}

We will be discussing a third family of groups as well.  These groups
are standard restricted wreath products of well known groups.  Wreath
products are discussed in detail in many places, but two standard
references are Peter Neumann's paper \cite{NeumannW} and the book
\cite{Meldrum} by J. D. P. Meldrum.  We now give a general description
of the standard restricted wreath product.

Let the symbol $C\wr T$ represent the standard restricted wreath
product of the groups $C$ and $T$.  In particular, $C\wr T\cong
B\rtimes T$, where $B\cong \bigoplus_{t\in T}C$.  We use the
coordinate system provided by the isomorphism of $C\wr T$ with
$B\rtimes T$ to describe elements of $C\wr T$.  In particular, if
$((a_t)_{t\in T},r)$, $((b_t)_{t\in T}, s)\in C\wr T$, then we
describe the product of these two elements by the following formula.
\[
((a_t)_{t\in T}, r)\cdot((b_t)_{t\in T},s) = ((c_t)_{t\in T}, rs)
\textrm{ where for }t\in T, \textrm{ we have }c_t = a_t
b_{(t\cdot r)}.
\]

In this notation, $C$ is called the bottom group, $B$ is called the
base group, and $T$ is called the top group.  We often think of the
base and top groups as being identified with the standard appropriate
subgroups in $B\rtimes T$, and we will sometimes refer to $B$ and $T$
as subgroups of $C\wr T$.
 
Given any positive integer $n$, let $ZW_n$ represent the group
$\Z\wr(\Z\times\Z\times\ldots\times\Z)$, where there are $n$
appearances of the factor $\Z$ in the direct product.

\subsection{Notation and theory associated with group actions and homeomorphisms of the circle}
All group actions in this paper will be right actions, that is, if $x$
is a point in a space $S$ acted on by a group $G$ and $a\in G$, then
we will denote by $xa$ the image of $x$ under the action of $a$.
Likewise, if $X\subset S$, we will denote by $Xa$ the set
$\left\{xa\,|\,x\in X\right\}$.  Consistent that notation, if $b$ is
also in $G$, we will use the following notational conventions, $a^b =
b^{-1}ab$ and $[a,b] = a^{-1}b^{-1}ab$.  Standard to this variety of
notation is the equation $Support(a^b) = (Support(a))b$, where the
support of an element of $G$ will be the set of points in $S$ moved by
the element.

There are three main results that we will use in our arguments, that
are relevant to the classic body of theory associated with groups of
homeomorphisms of the circle.  

The first result is the classical Ping-Pong Lemma, the proof of which
can be found in many places.  The earliest reference I know of is
\cite{pingpong}, a paper in German by Felix Klein.  The statement I
give follows that in \cite{FarbPingPong}.

\bl [Ping-Pong Lemma] Let $\Gamma$ be a group of permutations on a
set $X$, let $g_1$, $g_2$ be elements of $\Gamma$ of order at least
three.  If $X_1$ and $X_2$ are disjoint subsets of $X$ and for all
$n\neq 0$, $i \neq j$, $X_jg_i^n \subset X_i$, then $g_1$, $g_2$
freely generate the free group $F_2$.  
\el

While this lemma is clearly not restricted to groups acting on the
circle, it plays a heavy role in the theory of such groups,
particularly when the action of such a group is highly transitive.

 We will also make heavy use of Poincar\'e's concept of the rotation
number of a self-homeomorphism of $S^1$.

  Given a homeomorphism $f:S^1\to S^1$, let $F:\R \to \R$ represent a
lift of this homeomorphism via the standard covering projection
$exp:\R \to S^1$, so that $0 \leq F(0)<1$, and where we think of $S^1$
as a subset of the complex plane, and use $exp(t) = e^{2\pi i t}$.  We
then define
\[
\Rot{f} = lim_{n\to \infty}\frac{F^n(x)-x}{n}.
\]
Poincar\'e showed that this limit is well defined and it is independent
of the initial point $x\in[0,1)$ used.  Given an $f\in Homeo(S^1)$, we
call $\Rot{f}$ the \emph{rotation number of $f$}.  If $H$ is an
abelian subgroup of $Homeo(S^1)$, then the rotation number map
represents a homomorphism from $H$ to the additive group $\R \textrm{
mod }(1)$.

The second result we will need is the following lemma proved by Poincar\'e circa 1900.

\bl [Poincar\'e's Lemma] Suppose $f:S^1\to S^1$ is a homeomorphism with
rational rotation number $\frac{r}{s}$ (in lowest terms, use $r = 0$,
$s = 1$ if the rotation number is zero), then $f^s$ has a fixed point.
\el

The final result is a theorem by Denjoy, although we use an extension
given by M. Herman on page 76 of \cite{Herman}.

\bt [Denjoy's Theorem] Suppose $h:S^1\to S^1$ is either a
  piecewise-linear homeomorphism or a $C^1$ diffeomorphism whose
  derivative has bounded variation.  If the rotation number of $h$ is
  irrational, then $h$ is conjugate (by an element in $Homeo(S^1)$) to
  a rotation.

\et

\subsection{Statements of Results}
In \cite{bpgsc}, we see that $ZW_2$ does not embed in $\ploi$ (the
subgroup of $\pli{1}$ consisting only of orientation-preserving
elements), and therefore, that $ZW_2$ does not embed in R. Thompson's
group $F$, answering a question of Mark Sapir.  In a discussion of
that work, Ken Brown asked the author whether it could extended to get
a similar result for R. Thompson's group $T$.  We in fact show the
following theorem.

\bt
\label{nonembedding}
$ZW_2$ does not embed in $\pls{1}$.
\et

Since $T$ is often realized as a subgroup of $\pls{1}$, we immediately
obtain the following result.

\bc
$ZW_2$ does not embed in R. Thompson's group $T$.
\ec

Taking the above together with the following easy result (which we
leave for the reader's entertainment), we are lead to ask (in the next
sub-section) some natural follow-up questions.

\bt

Given any natural number $n>0$, $ZW_n$
embeds in $\pli{n}$, and also in $\pls{n}$.  

\et

\subsection{Two Questions}
The chief non-embedding result of \cite{bpgsc} shows
that $ZW_2$ fails to embed in $PL_o(I)$, which immediately implies
that $ZW_2$ fails to embed in $\pli{1}$.  The chief non-embedding
result here shows that $ZW_2$ fails to embed in $\pls{1}$.  The nature
of the proof of the non-embedding of $ZW_2$ in $\pli{1}$ leads the
author to suspect that these non-embedding results will generalize to
higher dimensions, but he has no idea as to how to approach any proof
of that suspicion.  Therefore, due to a lack of large quantities of
creditable evidence, we do not put forth any conjectures, but instead
we ask some questions.

\bq

Does there exist a natural number $n>1$ so that $ZW_n$ embeds in $\pli{n-1}$?

\eq

Note that given any piecewise-linear homeomorphism of the unit
$n$-cube, it is easy to build an induced piecewise linear
homeomorphism on the $n$-sphere by using the $n$-cube homeomorphism on
both halves of the sphere, which we recall is constructed as a doubled
$n$-cube.  In particular, we have the following remark.

\brk

For all $n\in\N$, the group $\pli{n}$ embeds in $\pls{n}$.

\erk

The remark still leaves open the possibility that there is an
$n\in\N$ so that $ZW_n$ embeds in $\pls{n-1}$, even if it does
not embed in $\pli{n-1}$.  Therefore, the following question is of
independent interest.

\bq

Does there exist $n\in\N$ so that $ZW_n$ embeds in $\pls{n-1}$?

\eq

Therefore, we have the following diagram, characterizing the
situation, where the symbols $PL(A^n)$ can be taken to represent
either $\pli{n}$ or $\pls{n}$.
\[
\xymatrix {
ZW_4\ar@{^{(}->}[rr]^{^{\vdots}}\ar@{^{(}->}[drr]^? & & PL(A^4)\\
ZW_3\ar@{^{(}->}[rr]\ar@{^{(}->}[drr]^{?} & & PL(A^3)\\
ZW_2\ar@{^{(}->}[rr]\ar@{^{(}->}[drr]^{!\exists} & & PL(A^2)\\
ZW_1\ar@{^{(}->}[rr] & & PL(A^1)
}
\]

\section{$ZW_2$ fails to embed in $\pls{1}$}
The proof which follows can be replaced by another proof, which relies
on the classification of the solvable subgroups of $\pls{1}$ which is
forthcoming in \cite{bkmss}.  However, the proof below is elementary, only
relying on the theory which was mentioned in the introduction.

\underline{Proof of Theorem \ref{nonembedding}}

Suppose there exists an injective homomorphism
$\phi:ZW_2\to \pls{1}$.  Recall that $ZW_2 = \Z\wr(\Z\times\Z)=\Z\wr(\Z^2)$.  Choose $a =
((0)_{*\in\Z^2},(1,0))$ and $b = ((0)_{*\in\Z^2},(0,1))$ as generators
of the top group $\Z^2$ in $ZW_2$, and let us use $c =
((c_*)_{*\in\Z^2},(0,0))$, (where $c_* = 1$ if $* = (0,0)$ and $c_* =
0$ otherwise) as our generator for an embedded copy of the bottom
group $\Z$ into the base group in $ZW_2$, located at index $(0,0)$ in
the base group.

Let $a\phi = \alpha$, $b\phi = \beta$, and $c\phi = \gamma$.  Our
argument will break into cases.

\vspace{.15 in} \underline{Case 1}: $\alpha$ or $\beta$ has a rational
rotation number.
\vspace{.1 in}

We will assume without meaningful loss of generality that $\alpha$ has
a rational rotation number.  If this number is not zero, then by
Poincar\'e's Lemma, there is an integer $n$ so that $\alpha^n$ has a
fixed point.  Further, as $\alpha$ has infinite order, we need not be
concerned that $\alpha^n$ is trivial.  Replace $a$ by $a^n$.  The new
$a$ and the old $b$ still generate a group isomorphic to $\Z^2$, and
$\langle a,b,c\rangle$ is still isomorphic with $ZW_2$.  Therefore, we
may assume that the rotation number of $\alpha$ is zero, and that
$\alpha$ is non-trivial, but still admits a fixed point $p$.

Let $r$ be any piecewise-linear homeomorphism of $S^1$ which sends $p$
to $1\in S^1$, so that $\alpha^r = r^{-1}\alpha r$ will be a
piecewise-linear homeomorphism which fixes $1$.  We can now replace
$\phi$ by a new embedding of $ZW_2$ into $PL(S^1)$ defined by the rule
$f\mapsto (\phi(f))^r$ for all $f\in ZW_2$.  From now on, when we
refer to $\phi$, we will mean the new version of $\phi$, and when we
refer to $\alpha$, $\beta$, or $\gamma$, we will be referring to the
images of $a$, $b$ , and $c$ under this new version of $\phi$.

Choose the representative lift $\tilde{\alpha}:\R\to\R$ of $\alpha$ so
that $0\tilde{\alpha} = 0$ (and therefore $1\tilde{\alpha} = 1$,
etc.).  We may also assume that we chose $p$ intelligently, so that
$\tilde{\alpha}$ has a component of support $\tilde{A}$ of the form
$(x,1)$ for some $x\in (0,1)$.  Let $A = exp(\tilde{A})$ be the image
of $\tilde{A}$ under the covering projection $exp:\R \to S^1$.  Now,
if $\beta$ has an irrational rotation number, by the extended version
of Denjoy's theorem, we can find a natural number $n$ so that
$1\beta^n= y\in A$.  The point $y$ represents a fixed point of the map
$\alpha^{(\beta^n)}$, where $\alpha$ has no fixed point.  In this
case, we see that $\alpha$ and $\beta$ fail to commute.  Hence $\beta$
must have a rational rotation number, and in fact, some power
$\beta^k$ of $\beta$ must fix the point $1$.  Replace $b$ by $b^k$,
$ZW_2$ by its isomorphic subgroup generated by $a$, the new $b$, and
$c$, and replace $\beta$ by the image of the new $b$ under $\phi$.

Now we have that $\alpha$ and $\beta$ fix the point $1\in S^1$,
although both elements have non-trivial components of support.  We
also understand that if $K$ is a component of support of $\alpha$,
then it is also a component of support of $\beta$, or else it is completely
disjoint from the support of $\beta$, and vice-versa.

Note that $[c,a] = c^{-1}a^{-1}ca$ is a non-trivial element in the
base group of $ZW_2$, and it has the property that $[c,a]$, $a$, and
$b$ generate a subgroup of $ZW_2$ which is isomorphic to $ZW_2$.  Let
$\theta$ be the image of $[c,a]$ under $\phi$.  The element $\theta$
is the product of two elements with fixed points, so if $\theta$ does
not have fixed points, then we must have that $\alpha$ and
$\alpha^{\gamma} = \phi(c^{-1}a^{-1}c)$ generate a group which
contains a free group on two generators (by Klein's classical
ping-pong argument (the earliest reference I know of is
\cite{pingpong}), where we use small neighborhoods of the ends of the
components of the support of the two elements to build the two
required sets for the ``Ping-pong'').  In particular, as $ZW_2$ does
not contain a free group on two generators, we see that $\theta$ must
admit fixed points.  We now replace $c$ by $[c,a]$, so that our new
$\gamma$ is the $\theta$ discussed in this paragraph.

If the set of elements $\Gamma = \left\{\alpha, \beta, \gamma\right\}$
admits a global fixed point, then it is easy to see that we can embed
the group $\langle \Gamma\rangle\cong ZW_2$ into $\pli{1}$, which we
know is impossible, by the chief non-embedding result of \cite{bpgsc}.
If $\Gamma$ does not admit a non-empty global fixed set, then we can find an
element $\theta$ in $\langle \alpha, \beta\rangle$ whose fixed set is
entirely contained in the support of $\gamma$, so that again a
ping-pong argument shows that $ZW_2$ must contain an embedded copy of
the free group on two generators.  

Hence, we have found contradictions in all cases, so that if $ZW_2$
embeds in $\pls{1}$, the generators $a$ and $b$ must be realized by
elements with irrational rotation numbers.

\vspace{.15 in}
\underline{Case 2}: $\gamma$ has a rational rotation number. 

We immediately replace $c$ by some non-trivial power of itself so that
$\gamma$ has a fixed point, and we let $p$ represent the end of some
component of support $D$ of $\gamma$.

Since $\alpha$
must have an irrational rotation number, the orbit of $p$ under the action
of $\alpha$ is dense in $S^1$, and so there is an integer power $k$ so
that $p\alpha^{k}\in D$.  But now the elements $\gamma$ and
$\gamma^{\alpha^k}$ cannot commute with each other, which contradicts
the fact that they are both in the image of the base group of $ZW_2$
under $\phi$, which group is abelian.  

Taking this contradiction together
with the results of the previous
case, we see that if $ZW_2$ embeds
in $\pls{1}$, we must have that all
of the generators $\alpha$, $\beta$,
$\gamma$ of the image have
irrational rotation number.
Further, observe that no non-trivial
element $\omega$ of the base group
can have image with rational
rotation number, or the group
$\langle \alpha, \beta,\omega\rangle
\cong ZW_2$ where the generator of
the new bottom group would have rational
rotation number.

\vspace{.1 in}

\vspace{.15 in} \underline{Case 3}: None of $\alpha$, $\beta$, or
$\gamma$ has a rational rotation number.
\vspace{.1 in}

Apply the extended version of Denjoy's theorem to find $\theta\in
Homeo(S^1)$ so that $\alpha^{\theta}$ is a pure rotation.  Replace
$\alpha$, $\beta$ and $\gamma$ by their images under conjugation by
$\theta$, so that the group they generate is now a subgroup of
$Homeo(S^1)$ (note that the image of $\gamma$ under this conjugation may not be
piecewise linear). 

If $\gamma$ is also a pure rotation, then $\alpha$ and
$\gamma$ would commute.  In particular, their is a power $n$ of
$\alpha$ so that the graphs of $\alpha^{-n}$ and $\gamma$ will intersect
each other on the torus (consider that the graphs of the powers of
$\alpha$ produce a dense subset of the foliation of the torus by the
graphs of the pure rotations).  In particular, the functions
$\alpha^{-n}\gamma^{-1}$ and $\alpha^n\gamma$ both have non-trivial
fixed sets.  

If the product $[\alpha^n,\gamma] =
\alpha^{-n}\gamma^{-1}\alpha^n\gamma$ has no fixed points, then we see
that the fixed set of $\alpha^{-n}\gamma^{-1}$ must be contained in the
support of $\alpha^n\gamma$, and vice versa, so that one can use a
ping-pong argument to show that $ZW_2$ contains a free group on two
generators.

If $[\alpha^n,\gamma]$ has some fixed points, then we have found a
non-trivial element in the base group of $ZW_2$ with rational rotation
number, which is impossible by the work in our previous case.

\qquad$\diamond$

\newpage
\bibliographystyle{amsplain}

\bibliography{ploiBib}

\end{document}